\documentclass[11pt]{article}
\usepackage{amsfonts}

\newtheorem{thm}{Theorem}

\newtheorem{lemma}[thm]{Lemma}
\newtheorem{prop}[thm]{Proposition}

\usepackage{epsfig}

\begin{document}

\def\G{{\Gamma}}
 \def\d{{\delta}}
 \def\ci{{\circ}}
 \def\e{{\epsilon}}
 \def\l{{\lambda}}
 \def\L{{\Lambda}}
 \def\m{{\mu}}
 \def\n{{\nu}}
 \def\o{{\omega}}
  \def\O{{\Omega}}
  \def\s{{\sigma}}
 \def\v{{\varphi}}
 \def\a{{\alpha}}
 \def\b{{\beta}}
 \def\p{{\partial}}
 \def\r{{\rho}}
 \def\ra{{\rightarrow}}
 \def\lra{{\longrightarrow}}
 \def\g{{\gamma}}
 \def\D{{\Delta}}
 \def\La{{\Leftarrow}}
 \def\Ra{{\Rightarrow}}
 \def\x{{\xi}}
 \def\c{{\mathbb C}}
 \def\z{{\mathbb Z}}
 \def\2{{\mathbb Z_2}}
 \def\q{{\mathbb Q}}
 \def\t{{\tau}}
 \def\u{{\upsilon}}
 \def\th{{\theta}}
 \def\la{{\leftarrow}}
 \def\lla{{\longleftarrow}}
 \def\da{{\downarrow}}
 \def\ua{{\uparrow}}
 \def\nwa{{\nwtarrow}}
 \def\swa{{\swarrow}}
 \def\nea{{\netarrow}}
 \def\sea{{\searrow}}
 \def\td{{\bigtriangledown}}
\def\hla{{\hookleftarrow}}
 \def\hra{{\hookrightarrow}}
 \def\sl2{{SL(2,\mathbb C)}}
 \def\ps{{PSL(2,\mathbb C)}}
 \def\qed{{\hfill$\diamondsuit$}}
 \def\pf{{\noindent{\bf Proof.\hspace{2mm}}}}
 \def\ni{{\noindent}}
 \def\sm{{{\mbox{\tiny $M$}}}}
 \def\sc{{{\mbox{\tiny $C$}}}}
 \def\ke{{\mbox{ker}(H_1(\p M;\2)\ra H_1(M;\2))}}
 \def\et{{\mbox{\hspace{1.5mm}}}}
\def\sk{{{\mbox{\tiny $K$}}}}
 \def\skp{{{\mbox{\tiny $K_+$}}}}
 \def\skm{{{\mbox{\tiny $K_-$}}}}
  \def\sb{{{\mbox{\tiny $\b$}}}}
 \def\sd{{{\mbox{\tiny $\d$}}}}
\def\sv{{{\mbox{\tiny $V$}}}}
\def\sl{{{\mbox{\tiny $\l$}}}}

\begin{center}
{\bf  Positive Knots And Knots With Braid Index Three Have Property-P }
\end{center}

\begin{center}
W. Menasco   and X. Zhang\footnote{{Partially supported by NSF
grant DMS 9971561}}
\end{center}

{\small
\begin{center} Department of Mathematics, SUNY at Buffalo, Buffalo, NY
14260-2900

 menasco@math.buffalo.edu \hspace{4mm}and \hspace{4mm}
xinzhang@math.buffalo.edu
\end{center}

{\bf Abstract:} We prove that positive knots and knots with
braid index three in the 3-sphere satisfy the Property P conjecture. }

 \vspace{5mm}

\begin{center}
{\bf 1. Introduction}
\end{center}

Let $K$ be a knot in the 3-sphere $S^3$ and $M = M_\sk$ the complement
of an open regular neighborhood of $K$ in $S^3$.  As usual, the set of
slopes on the torus $\p M$ (i.e. the set of isotopy classes of
unoriented essential simple loops on $\p M$) is parameterized by
$$\{m/n \,; \,\, m,n\in \z, n>0, (m,n)=1\}\cup \{1/0\},$$ so
that $1/0$ is the meridian slope of $K$ and $0/1$ is the longitude slope of $K$.
The manifold
obtained by Dehn surgery on $S^3$ along the knot $K$ (equivalently,
Dehn filling on $M$ along the torus $\p M$) with slope $m/n$, is
denoted by $K(m/n)$ or $M(m/n)$.
Of course $K(1/0)=S^3$, and thus the surgery  with
the slope $1/0$ is called the trivial surgery.
The celebrated  {\it Property P} conjecture, introduced by
Bing and Matin in 1971 [BMa],  states that
every nontrivial  knot $K$ in $S^3$ has Property P, i.e.
every nontrivial surgery on $S^3$ along $K$ produces a non-simply
connected manifold.
For convenience we say that  a class of knots in $S^3$ have
Property P
if every nontrivial knot in this class has Property P.
The following classes of knots were known to have Property P:
torus knots [H],  symmetric knots [CGLS] (the part
for strongly invertible knots
was proved in [BS]),
satellite knots [G1], arborescent knots [W], alternating knots [DR],
and small knots with no non-integral boundary slopes [D].
 For a simple homological reason,
to prove the conjecture for a knot $K$
one only needs to consider the  surgeries of $K$ with slopes $1/n$, $n\ne 0$.
A remarkable progress on the conjecture was made in [CGLS];
it was proved there that for a nontrivial knot, only one of $K(1)$ or $K(-1)$  could possibly
be a simply connected manifold.
Another remarkable result was given in  [GL], which told us that
 if the Property P conjecture is false, then  the Poincare conjecture is false.
For some earlier progresses on the conjecture, see [K, Problem 1.15]
for a summary. In this paper we prove

\begin{thm}\label{thm1}
Positive (or negative) knots in $S^3$ satisfy the Property P conjecture.
\end{thm}

Recall that a knot is positive if it can be represented as the closure
of a positive $n$-braid for some $n$, i.e. a braid which
involves the standard elementary braid generators
$\sigma_1,\cdots,\sigma_{n-1}$ (Figure 1 gives $\s_1,\s_2$ when $n=3$)
but not their inverses.
A negative knot is similarly defined, which is actually
the mirror image of a positive knot.

\begin{thm}\label{thm2}
Knots in $S^3$ with braid index three satisfy the Property P conjecture.
\end{thm}

Theorem 1 is a quick application of the Casson invariant. We refer
to [AM] for the definition and basic properties of the Casson
invariant. The Casson invariant is an integer valued topological
invariant defined for homology 3-spheres and for knots in homology
3-spheres. If a homology 3-sphere has non-zero Casson invariant,
then the manifold has an irreducible representation from its
fundamental group to the group $SU(2)$, which implies in
particular that the manifold
 is non-simply connected.
 For a
knot $K$ in $S^3$, if $\D_\sk(t)$ denotes the normalized Alexander
polynomial of $K$, i.e. satisfying $\D_\sk(1)=1$ and
$\D_\sk(t^{-1})=\D_\sk(t)$, then the Casson invariant $C_\sk$ of $K$ is
equal to the integer $\frac12\D_\sk''(1)$ [AM].
Further the Casson invariant of the  manifold  $K(1/n)$
is  equal to $n C_\sk$. So  $C_\sk\ne 0$
implies in particular  that the knot $K$ has Property $P$. Note that the Conway
polynomial of a knot in $S^3$ is a single variable  polynomial
$\td_\sk(x)\in \z[x]$ with only even powers and
$\td_\sk(t^{1/2}-t^{-1/2})=\D_\sk(t)$, from which  one can easily deduce
that the coefficient of $x^2$ in $\td_\sk(x)$ is equal to the
Casson invariant of $K$. For a nontrivial positive knot $K$ in
$S^3$, it has been proved in [V] that the coefficient of $x^2$ in
$\td_\sk(x)$ is a positive integer. Hence such knot has  Property P.
A negative knot is just a mirror image of a positive knot, and
obviously a knot has Property P if and only if its
 mirror image does.
Theorem \ref{thm1} now follows.

The rest of the paper is devoted to the proof of Theorem
\ref{thm2}. The main tools we shall use are the Casson invariant and
essential laminations.
We refer to [GO] for the definition and basic properties of
an essential lamination.
 In section 2 we give an outline of the proof of
Theorem \ref{thm2}. Actually the proof of Theorem \ref{thm2} is reduced
there to that of two propositions, Proposition 3 and Proposition 4.
These two propositions will then be proven  in section 3 and
section 4 respectively.

\vspace{2mm}
\begin{center}
{\bf 2. Proof of Theorem \ref{thm2}}
\end{center}

Recall that the 3-braid group, $B_3$,  has the following well known
 Artin presentation: $$B_3=<\s_1,\s_2\;|\;\s_1\s_2\s_1=\s_2\s_1\s_2>$$
where $\s_1$ and $\s_2$ are elementary 3-braids as shown in Figure 1.
\begin{figure}[!ht] {\epsfxsize=1.5in
\centerline{\epsfbox{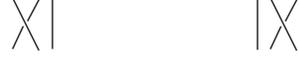}}} \caption{$\s_1$ (the left figure)
and  $\s_2$ (the right figure)}
\end{figure}
If we let $a_1=\s_1, a_2=\s_2,
a_3=\s_1^{-1}\s_2\s_1=\s_2\s_1\s_2^{-1}$, then $B_3$ also has the
following presentation in generators $a_1,a_2$ and $a_3$ (see
[X]): $$B_3=<a_1,a_2,a_3\;|\; a_2a_1=a_3a_2=a_1a_3>.$$ In this paper we shall
always express a 3-braid as
 a word  $w(a_1,a_2,a_3)$ in letters $a_1,a_2,a_3$. Such
a word is called {\it positive}  if the power of every letter
in the word is positive. A positive
 word $w=a_{\t_1}\cdots  \a_{\t_k}$ is said to be in
{\it  non-decreasing order} ($ND$-order) if the array of its
subscripts $(\t_1, ... ,\t_k)$ satisfies $$\t_{j+1}=\t_j\;\;
\mbox{or} \;\; \t_{j+1}=\t_j+1\;(\mbox{mod $3$ if $\t_j+1=4$})\;\;\mbox{for
$j=1,...,k-1.$}$$
One can define {\it negative} word and {\it
 non-increasing order} ($NI$-order) similarly. Let $P$ be the set
of positive words in $ND$-order, let $N$ be the set of negative
words in $NI$-order, and let $\a=a_2a_1$.
 It is proven in [X]
that for any 3-braid, there is a representative in its conjugacy
class that is a shortest word in $a_1, a_2, a_3$ and is
 of the form

\noindent
(i) a product of $\a^k$ and a word (maybe empty) in $P$
for some non-negative integer $k$; or

\noindent (ii)  a product of $\a^{k}$ and a word (maybe empty) in
$N$ for some non-positive integer $k$; or

\noindent(iii)  a product of   a word in $P$ and a word in $N$.

\noindent
where the meaning  of {\it the shortest} is
that the length, i.e. the number of letters of the representative
is minimal among all  representatives
 in the conjugacy class of the braid.
Such representative of a 3-braid
is said to be in  {\it normal form}. We shall  only need to show
 that if $K$ is a nontrivial knot in $S^3$ which is the  closure of a
 3-braid in normal form (i) or (ii) or (iii), then it has Property P.

  Recall that a word $w(a_1,a_2,a_3)$ is called {\it freely
reduced} if no adjacent letters are inverse to each other, and is
called {\it cyclically reduced} if it is freely reduced and the
first letter and the last letter of the word are not inverse to
each other. Given a word $w(a_1,a_2,a_3)$, one can combine all
adjacent letters of the same subscript into a single power of the
letter, called a {\it syllabus} of the word in that subscript. A
word $w$ is called {\it syllabus reduced} if it is expressed as a
word in terms of syllabuses as
$$w=a^{\mbox{\tiny{$m_1$}}}_{\mbox{\tiny{$\t_1$}}}
a^{\mbox{\tiny{$m_2$}}}_{\mbox{\tiny{$\t_2$}}}\cdots
a^{\mbox{\tiny{$m_k$}}}_{\mbox{\tiny{$\t_k$}}} $$ such that
$a_{\mbox{\tiny{$\t_j$}}}\ne a_{\mbox{\tiny{$\t_{j+1}$}}}$ for
$j=1,...,k-1$.  A word is called {\it cyclically syllabus reduced}
if it is syllabus reduced and its first and last syllabuses are
 in different subscripts.

Let $P^*$ denote the set of all positive words in $a_1,a_2,a_3$
such that between any two syllabuses in $a_3$ both $a_1$ and $a_2$
occurs. Obviously any 3-braid of norm form (i) is contained in
$P^*$.
Let $\b$ be a 3-braid in $P^*$. Suppose that $a_3^k$ is a
syllabus in $\b$ which is proceeded immediately by $a_1$.
Then one can eliminate the syllabus $a_3^k$ with the equality
$a_1a_3^k=a_2^ka_1$ to get an isotopic braid which is still in
$P^*$  but with one less  number of syllabuses in $a_3$.
Similarly if a syllabus $a_3^k$ is followed immediately by
$a_2$, then one can eliminate the syllabus $a_3^k$ with the
equality $a_3^ka_2=a_2a_1^k$ to get an isotopic braid which is
still in $P^*$  but with one less number of syllabuses in $a_3$.
We shall call this process {\it index-3 reduction}. So for any
given $\b\in P^*$, we can find, after a finitely many times of
index-3 reduction, an equivalent braid representative $\b'$ in
$P^*$ for $\b$ such that
 every syllabus in $a_3$ occurring in $\b'$ can only
possibly be
proceeded immediately  by $a_2$ and likewise can only possibly
be followed immediately by $a_1$.
We call a  word $\b$ in $P^*$ {\it  index-3 reduced} if
  every syllabus in $a_3$ occurring in $\b$ is neither
proceeded immediately  by $a_1$  nor followed immediately by
$a_2$.

Let $P^a$ denote the set of 3-braids of the form
 $\b=a_i^{-q}\d$, where $q=0$ or $1$
 and $\d\in P^*$ is a non-empty word, such that
 $\b$  is cyclically reduced and $\d$ is index-3 reduced.
Obviously $P$ is contained in $P^a$.

\begin{prop}\label{casson}
Suppose that $K$ is a  knot in $S^3$ which is the closure of a
3-braid $\b=a_i^{-q}\d$ in $P^a$ such that $\d$ contains at least
four syllabuses but $\b$ is not one of the words in the set $E=$\{
$a_1^{-1}a_2a_3^2a_1a_2$, $a_1^{-1}a_3^2a_1a_2a_3$,
 $a_1^{-1}a_3a_1a_2^2a_3$, $a_1^{-1}a_2a_3a_1a_2^2$,
  $a_2^{-1}a_3a_1a_2a_3^2$,  $a_2^{-1}a_3a_1^2a_2a_3$,
  $a_2^{-1}a_1a_2a_3^2a_1$, $a_2^{-1}a_1^2a_2a_3a_1$,
  $a_3^{-1}a_1a_2a_3a_1^2$, $a_3^{-1}a_1a_2^2a_3a_1$,
  $a_3^{-1}a_2a_3a_1^2a_2$, $a_3^{-1}a_2^2a_3a_1a_2$\}.
  Then $K$ has positive Casson
invariant and thus has Property P.
 \end{prop}

Later on we shall refer  the set $E$ given
in Proposition \ref{casson} as the {\it excluded} set.

\begin{prop}\label{lamin}
Suppose that $K$ is a  knot in $S^3$ which is the closure of a
3-braid $\b=\d\eta$ which is  in normal
form (iii), i.e. $\d\in P$ and $\eta\in N$.  Suppose that
either
\newline
(1) each of $\d$ and $\eta$ has a syllabus of power larger than one,   or
\newline
(2) each of  $\d$ and
$\eta$ contains at least two syllabuses, or
\newline
(3) one of $\d$ and $\eta$  contains at least four
syllabuses and the other has length at least two.
\newline
Then each of $K(1)$ and $K(-1)$ is a manifold
which contains an essential lamination.
\end{prop}

If a closed 3-manifold has an essential
lamination, then its universal cover is $\mathbb R^3$ [GO] and
thus in particular
the manifold cannot be simply connected.
Hence any knot as given in Proposition \ref{lamin}
has Property P by [CGLS].

Given Propositions \ref{casson} and \ref{lamin}, we can finish the
proof of Theorem \ref{thm2} as follows. 
For a braid $\b$, we use $\hat \b$ to denote the closure of $\b$.
Let $K\subset S^3$ be a
nontrivial knot with index $3$. Let $\b$ be a 3-braid
in normal form such that $\hat \b=K$.

First we consider the case that  $\b$ is in normal form (i), i.e.
$\b=\a^k\d$ with $\d\in P$ and $k\geq 0$. If $\b$ contains at
least four syllabuses and belongs to $P^a$, then $K=\hat\b$  has positive
Casson invariant by Proposition \ref{casson}. So the knot $K$ has
Property P in this case. If $\b$ has less than four syllabuses, then
 up to conjugation in $B_3$,  $\b=a_2a_1a_3^i$ or
 $\b=a_1^ia_2^ja_3^m$ for $i,j,m\geq 0$.
It is easy to see that in this case $\hat\b$ is an arborescent
knot and thus by [W], $K=\hat\b$ has Property P. For instance if
$\b=a_1^ia_2^ja_3^m$, then $\hat\b$ is as
shown in Figure 2 which shows in fact that $\hat\b$ is a
Montesinos knot. \begin{figure}[!ht] {\epsfxsize=5in
\centerline{\epsfbox{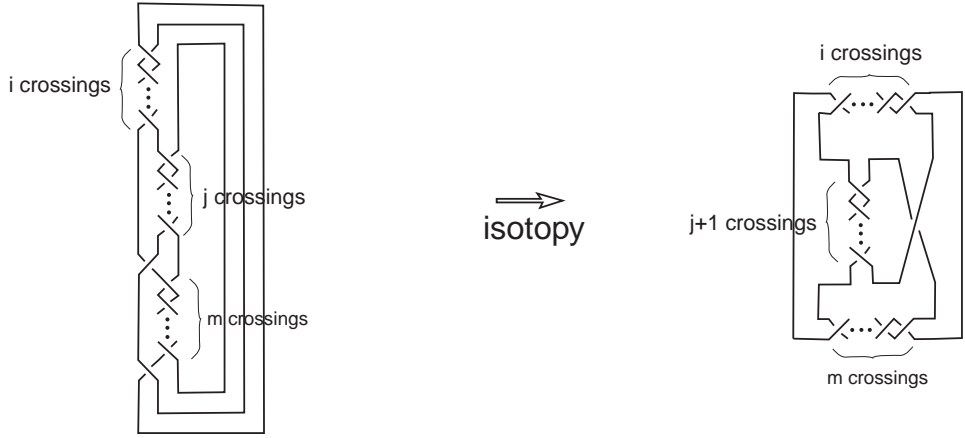}}} \caption{the closure of
$a_1^ia_2^ja_3^m$ is a Montesinos link}
\end{figure}
 So we may assume
that $\b$ has at least four syllabuses but is not in $P^a$. This
implies that in $\b=\a^k\d$, we have $k>0$ and  $\d$ starts  with
a syllabus in $a_3$. Performing index-3 reduction on $\b$ once, we
get an equivalent 3-braid $\b'\in P^*$ which is index-3 reduced,
i.e. $\b'\in P^a$. If $\b'$  has less than four syllabuses, then
again $K=\hat\b=\hat\b'$ is an arborescent knot and thus has
Property P. If $\b'$ contains at least four syllabuses, we may
apply Proposition \ref{casson} to get Property P for the knot.

If $\b$ is of  normal form
(ii), then the mirror image of $\b$ is a braid of normal form (i)
and thus the knot $K=\hat\b$ has Property P.

\begin{figure}[!ht] {\epsfxsize=5in \centerline{\epsfbox{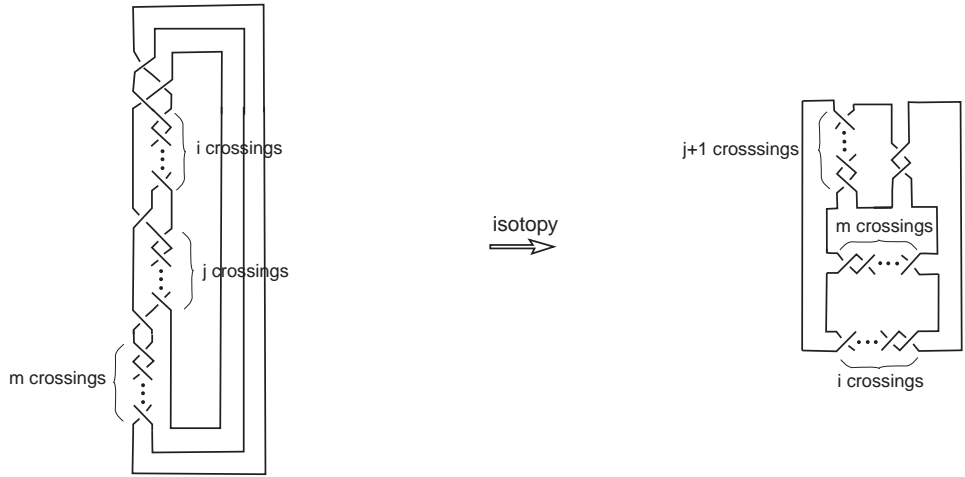}}}
\caption{the closure of $a_3^{-1}a_2^ia_3^ja_1^m$ is an
arborescent link}
\end{figure}

Suppose finally that $\b=\d\eta $ is a braid  of normal form (iii)
(recall that $\d\in P$ and $\eta\in N$). By Proposition
\ref{lamin}, we may  assume that each of the conditions (1)-(3) in
Proposition \ref{lamin} does not hold for $\b=\d\eta$. Then
  $\b$ is a word of the form $\d a_i^{-1}$ or $a_i\eta$,
 or $\b=\d a_i^{-k}$ where $k>1$ and $\d$  contains at most three
  syllabuses each having  power $1$,
 or $\b=a_1^{k}\eta$ where $k>1$ and $\eta$  contains at most three
  syllabuses each having  power $-1$.
If $\b$ contains at most four syllabuses,  then
 $K=\hat \b$ is an arborescent knot and thus has Property P.
(Figure 3 shows this for the case that
$\b=a_3^{-1}a_2^ia_3^ja_1^m$. Other cases can be treated similarly).
So we may only consider the cases
when
 $\b=\d a_i^{-1}$ or $\b=a_i\eta$,
  each  containing at least five syllabuses.
If $\b=\d a_i^{-1}$,  then $\b$ is conjugate to $\b'=a_i^{-1}\d$
which is in $P^a$. Hence if $\b'$ does not belong to the excluded
set $E=$\{ $a_1^{-1}a_2a_3^2a_1a_2$, $a_1^{-1}a_3^2a_1a_2a_3$,
 $a_1^{-1}a_3a_1a_2^2a_3$,$a_1^{-1}a_2a_3a_1a_2^2$,
  $a_2^{-1}a_3a_1a_2a_3^2$,  $a_2^{-1}a_3a_1^2a_2a_3$,
  $a_2^{-1}a_1a_2a_3^2a_1$, $a_2^{-1}a_1^2a_2a_3a_1$,
  $a_3^{-1}a_1a_2a_3a_1^2$,$a_3^{-1}a_1a_2^2a_3a_1$,
  $a_3^{-1}a_2a_3a_1^2a_2$, $a_3^{-1}a_2^2a_3a_1a_2$ \}, then
$K=\hat\b'$ has positive Casson invariant by Proposition
\ref{casson} and thus has Property P. If $\b'$ is in
the excluded set $E$,  then
 $K=\hat \b'$ is  an arborescent knot
 and thus has Property P [W].
(Figure 4 illustrates this for the case $\b=a_1^{-1}a_3^2a_1a_2a_3$.
Other cases can be checked similarly).
 Finally, if $\b=a_i\eta$,  then its mirror image is a braid
 in $P^a$, which is a case we have just discussed.
 This completes the
 proof of Theorem \ref{thm2}.

Propositions \ref{casson} and \ref{lamin} will be proved in
subsequent  two sections which constitutes the rest of the paper.

\begin{figure}[!ht] {\epsfxsize=3in
\centerline{\epsfbox{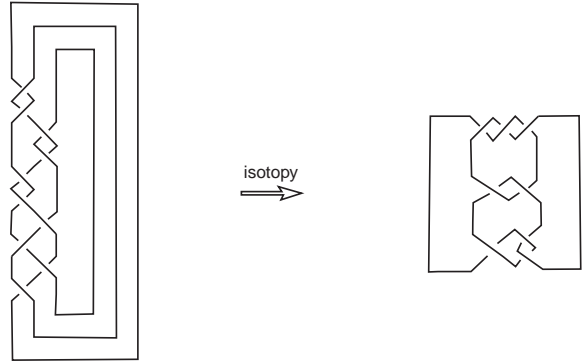}}} \caption{the closure of the braid
$a_1^{-1}a_3^2a_1a_2a_3$ is a Montesinos knot}
\end{figure}

\begin{center}
 {\bf 3. Proof of Proposition \ref{casson}}
\end{center}

We retain all definitions and notations established earlier. For a
3-braid $\b$, we use $n_\sb$ to denote the number of syllabuses in
$a_3$ occurring in $\b$ and use $s_\sb$ to denote the number of
syllabuses of $\b$. Obviously if $\b'$ is the braid obtained from
$\b\in P^*$ after some non-trivial index-3 reduction, then $\hat
\b'=\hat\b$ but $n_{\sb'}<n_\sb$. For a syllabus
$a_3^k=a_1^{-1}a_2^ka_1$, we shall always assume its plane
projection corresponds naturally to $a^{-1}_1a_2^ka_1$ as shown in
Figure 6 (a).   Hence, every 3-braid in letters $a_1,a_2,a_3$ has
its {\it canonical} plane projection; namely in the projection
plane  we place vertically (from top to bottom) and successively
the projections of letters   occurring in the braid, corresponding
to their natural  order from left to right. Whenever we need
consider a plane projection of a 3-braid, the canonical one is
always assumed unless  specifically indicated otherwise.

The basic tool we are going to use to prove the proposition
is the following crossing change formula $(*)$  of
the Casson invariant.
 If $K_+$, $K_-$ are oriented knots and $L_0$ an
 oriented link with two components in $S^3$ such that
 they have identical
plane projection except at one crossing they differ as shown in
Figure 5, then the Casson invariant of $K_+$ and $K_-$ satisfy the
following relation:
$$C_\skp-C_\skm=lk(L_0)\;\;\;\;\;\;\;\;\;\;\;\;\;\;\;\;\;\;\;
\;\;\;\;\;\;\;\;(*)$$ where $lk(L_0)$ is the linking number of
$L_0$. This formula can be found on page 141 of [AM] (note that
there was a print error there, $K_+$ and $K_-$ in Figure 36 and
Figure 37 of [AM] should be exchanged).
 \begin{figure}[!ht]
{\epsfxsize=3in \centerline{\epsfbox{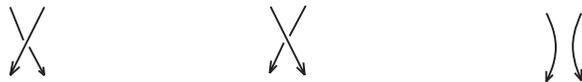}}} \caption{$K_+$ (the
left figure), $K_-$ (the middle figure)  and  $L_0$ (the right
figure)}
\end{figure}
The idea of proof of the proposition is repeatedly applying the
formula $(*)$ to a 3-braid of the type as given in the proposition
to reduce the complexity of the braid  and inductively prove its
positivity.
 To do so, we first need to  estimate the linking number
 of a two-component link which is the closure of
  a 3-braid  of relevant type.
 Given an oriented link $L$ of two components
 $L_1$ and $L_2$, we shall
calculate  the linking number of $L$ as follows (cf. [Rn]).
Take a plane projection of $L$.
 A crossing of $L$ as shown on the
 left of  Figure 5
 has  positive sign $1$ and the crossing in the middle of the figure
 has  negative sign $-1$.
The linking number of $L$ is the algebraic sum of the
crossing signs at those  crossings
of $L$  where $L_1$  goes under
$L_2$.

For a 3-braid $\b$, we shall always orient each component of $\hat\b$
in such a way  that the induced orientation on each strand of $\b$
in its canonical projection is pointing downward in the  projection plane.

Suppose that  a two component link $L=L_1\cup L_2$ is the closure
of a 3-braid $\b$. At a crossing corresponding to a letter
$a_1$ or $a_2$ or $a_1^{-1}$ or $a_2^{-1}$ appeared in $\b$, if
the under strand is from $L_1$ and the upper strand is from $L_2$,
then this crossing contributes negative one to the linking number
of $L=\hat\b$ when the crossing is corresponding to $a_1$ or $a_2$
and positive one when the crossing is corresponding to $a_1^{-1}$
and $a_2^{-1}$; and if the under strand is from $L_2$ or the upper
strand is from $L_1$, then this crossing contributes zero to the
linking number of $L=\hat\b$. The linking number of $L=\hat \b$
is the sum of the
contributions from all the crossings of $\b$. If $\b'$ is a
portion of $\b$, we use $l(\b')$ to denote the total contribution
to the linking number of $L$ coming from all the crossings of
$\b'$. In particular $l(\b)=lk(L)$. Also if we decompose $\b$ into portions
$\b=\b_1\b_2\cdots  \b_k$,
then $l(\b)=l(\b_1)+l(\b_2)+\cdots  +l(\b_k)$.

\begin{lemma}\label{lk0}
Let $\d\in P^*$ be a 3-braid  satisfying: (1) $n_\sd=0$; (2) $\hat\d$
is a link of two components; and (3) both $a_1$ and $a_2$ appear in
$\d$. Then $lk(\hat\d)<0$.
\end{lemma}

\pf  Recall that $n_\sd=0$ means that  $\d$ is a word in letters
$a_1$ and $a_2$ only. Also, $\d\in P^*$ is a positive word. The
conclusion  of the lemma is now obvious. \qed

\begin{lemma}\label{lk}
Let  $\b=a_i^{-1}\d$ be an element in $P^a$ satisfying: (1)
$n_\sd=0$; (2) $\hat\b$ is a link of two components; and (3) $\d$
contains at least four syllabuses. Then $lk(\hat \b)<0$.
\end{lemma}

\pf Again $\d$ is a positive word without letter $a_3$.
 If $i=1$,  then $\d$ must start with a syllabus in
$a_2$ and also ends with a syllabus in $a_2$ since $\b$ is
cyclically reduced. Since $\d$ contains at least four syllabuses,
it follows that $\b=a_1^{-1}a_2^ja_1^ka_2^ma_1^na_2^p\cdots  $ for
some $j,k,m,n,p>0$. One can easily verify that $lk(\hat \b)\leq
l(a_1^{-1}a_2^ja_1^ka_2^ma_1^na_2^p)<0$. Similarly one can treat
the $i=2$ case. Consider now the case $i=3$. If $\d$ starts with
$a_2$,  then $\b=a_3^{-1}a_2^ja_1^ka_2^ma_1^n\cdots  $ for some
$j,k,m,n>0$, which is conjugate to
$a_1^{-1}a_2^{j-1}a_1^ka_2^ma_1^n\cdots a_2$. So if $j>1$,  then
we are back to the case $i=1$ (all required conditions remain
valid). If $j=1$,  then $\b$ is conjugate to
$\b'=a_1^{k-1}a_2^ma_1^n\cdots  a_2$. Hence by Lemma \ref{lk0} we
see that $lk(\hat\b) =lk(\hat\b')<0$. Similarly one can treat the
case when $\d$ starts with $a_1$. \qed

\begin{lemma}\label{lk1}
Let  $L=L_1\cup L_2$ be a link of two components in $S^3$ which is the closure
of a braid $\b=a_i^{-q}\d$ in $P^a$ such that $\d$ contains at
least four syllabuses. Then the linking number $lk(L)$ of $L$ is
non-positive.
\end{lemma}

\pf
We may assume that $\b$ has been chosen in its conjugacy class,
subject to satisfying all conditions of the lemma,
to have the minimal number $n_\sd$.
 By Lemmas \ref{lk0} and \ref{lk}, we may assume that
$n_\sd>0$.

 Given the canonical plane projection of
a 3-braid $\eta$, we shall always call the strand of $\eta$ which starts
at the top left corner {\it the strand 1} of $\eta$,
call the strand which starts
at the top middle place {\it the strand 2} of $\eta$,
and call the strand which starts
at the top right corner {\it the strand 3} of $\eta$.

Let $a_3^k=a_1^{-1}a_2^ka_1$ be a syllabus in $a_3$ occurring in
$\d$. Its crossings are shown in Figure 6 (a).
\begin{figure}[!ht] {\epsfxsize=4in \centerline{\epsfbox{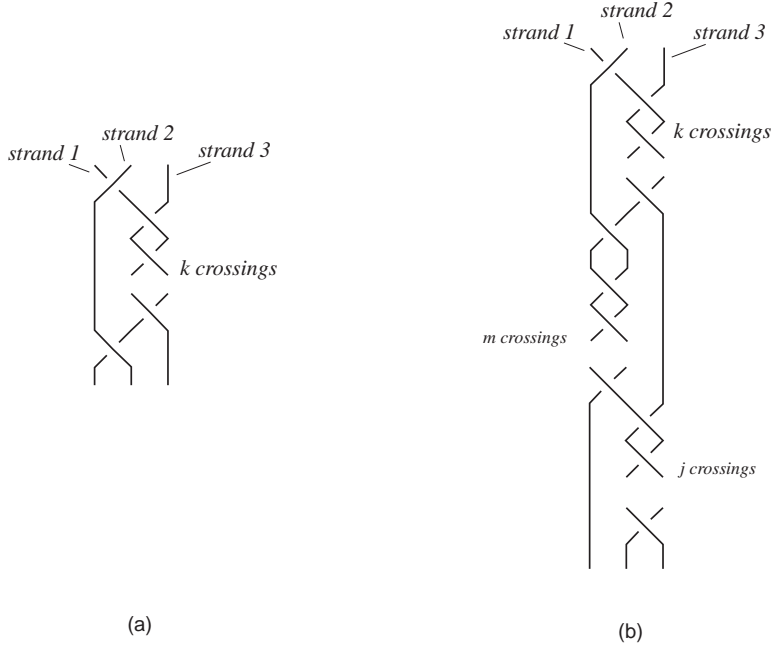}}}
\caption{(a) $a_3^k$\hspace{10mm} (b) $\d_0=a_3^ka_1^ma_2^j$}
\end{figure}
Obviously the total contribution of the syllabus to
the linking number is at most one; i.e. $l(a_3^k)\leq 1$.
More precisely, one can easily verify that the following claim
holds.

\noindent {\bf Claim 1}. If $l(a_3^k)=1$,  then the strand 1 of
$a_3^k$ (Figure 6 (a)) is from $L_1$, the strand 2 of $a_3^k$ is
from $L_2$,
 $k=1$, and the strand 3  is from $L_2$.

\noindent
{\bf Claim 2}. If $\d_0=a_3^ka_1^ma_2^j$ is a portion of $\d$ for some
$k, m, j >0$,
then $l(\d_0)\leq 0$.

Since $l(a_3^ka_1^ma_2^j)= l(a_3^k)+l(a_1^ma_2^j)$, the claim is
obviously true if $l(a_3)\leq 0$. So suppose that $l(a_3^k)=1$.
Then by Claim 1, we have that the strand 1 of $\d_0$ (Figure 6
(b)) is from $L_1$, $k=1$  and the strands 2 and 3 of $\d_0$ are
from $L_2$.
 So there
will be a negative one contribution to the linking number at
the first crossing of the
syllabus in $a_2$. This proves the claim.

We also need to consider under what conditions we have
$l(\d_0)=0$. Consider  Figure 6 (b).

\noindent Case (A1). If the strand 1 of $\d_0$ is from $L_1$ and the strands 2
and 3 of $\d_0$ are from $L_2$,  then $l(\d_0)=0$ happens exactly
in the following conditions:
 $k=1$ and $j\leq 2$;

\noindent Case (A2). If the strand 2 of $\d_0$ is from $L_1$ and
the strands 1 and 3 of $\d_0$  are from $L_2$,  then $l(\d_0)=0$
never happen.

\noindent Case (A3). If the strand 3 of $\d_0$ is from $L_1$ and
the strands 1 and 2 of $\d_0$ are from $L_2$,  then $l(\d_0)=0$
never happen.

\noindent Case (A4). If the strand 1 of $\d_0$ is from $L_2$ and
the strands 2 and 3 of $\d_0$ are from $L_1$,  then $l(\d_0)=0$
never happen.

\noindent Case (A5). If the strand 2 of $\d_0$ is from $L_2$ and the strands 1
and 3 of $\d_0$ are from $L_1$,  then $l(\d_0)=0$ happens exactly
in the following condition:
  $m=1$.

\noindent Case (A6). If the strand 3 of $\d_0$ is from $L_2$ and
the strands 1 and 2 of $\d_0$ are from $L_1$,  then $l(\d_0)=0$
never happen.

 Suppose that $\d_1$ is a
portion of the braid $\d$  which starts with a syllabus
in $a_3$ and is followed by a positive word
which starts
with $a_1$ and contains both  $a_1$ and $a_2$, but not $a_3$.
It follows from Claim 2 that

\noindent {\bf Claim 3.} $l(\d_1)\leq 0$.

We also need to know exactly when $l(\d_1)=0$. By the  discussion
following Claim 2, we  have that the portion of $\d_1$ which
consists of the first three syllabuses of $\d_1$ must be  as in
Cases A1 or  A5. One can then  verify that $l(\d_1)=0$ happens
exactly in one of the following situations:

\noindent
(B1) $\d_1=a_3a_1^ma_2^j$, for some $m>0, 0<j\leq 2$,
and the strand 1 of $\d_1$ is from $L_1$ and the strands 2 and 3
 from $L_2$; or

\noindent
(B2) $\d_1=a_3a_1^ma_2^2a_1^p$, for some $m,p>0$,and the
strand 1 of $\d_1$ is from $L_1$ and the strands 2 and 3 from $L_2$; or

\noindent
(B3) $\d_1=a_3^ka_1a_2^j$, for some $k,j>0$, and
the strands 1 and 3 of $\d_1$  are  from $L_1$ and the strand 2  from $L_2$.

Suppose that $\d_2$ is a portion of the braid $\d$ between and including
two syllabuses in $a_3$.

\noindent {\bf Claim 4}. $l(\d_2)\leq 0$.

We have $\d_2=\d_1a_3^r$ for some $r>0$ where $\d_1$ is a braid
as in Claim 3 and $\d_1$ must also  end with a syllabus in $a_2$ since $\d$
is index-3 reduced.
Suppose otherwise,  that is
$l(\d_2)>0$. We will get a contradiction.
By Claim 3 and Claim 1, we may assume
 that $l(\d_1)=0$ and $l(a_3^r)=1$.
So $\d_1$ is as in case (B1) or (B3)
and $a_3^r$ is as in Claim 1.
One can easily check that the link component assignment
to the strands of $\d_1$ and the link component assignment
to the strands of $a_3^r$ never match.

Similarly one can show

\noindent {\bf Claim 5}.  Suppose that  $\d_3$ is a portion of
$\d$ satisfying (1)  it ends with a syllabus in $a_3$; (2) this
ending syllabus is proceeded
by a positive word in $a_1$ and $a_2$ and the word contains
at least two syllabuses. Then $l(\d_3)\leq 0$ and $l(\d_3)=0$ happens
exactly in one of the following situations:

\noindent (C1) $\d_3=a_1a_2^ma_3$, for some $m>0$, and the strand
2 of $\d_3$  is from $L_1$ and the strands 1 and 3
 from $L_2$; or

\noindent
(C2) $\d_3=a_1^2a_2^ma_3$, for some $m>0$,
and the strand 1 of $\d_3$ is from $L_1$ and the strands 2 and 3
 from $L_2$; or

\noindent
(C3) $\d_3=a_2^pa_1^2a_2^ma_3$, for some $m,p>0$, and the
strand 1 of $\d_3$ is from $L_1$ and the strands 2 and 3 from $L_2$; or

\noindent
(C4) $\d_3=a_1^ja_2a_3^k$, for some $k,j>0$, and
the strands 1 and 2 of $\d_3$  are  from $L_1$ and the strands 3
from $L_2$.

\noindent {\bf Claim 6}.  Suppose that  $\d_4=\o\d_1$ is a portion
of $\d$ such that $\d_1$ is a braid as in Claim 3 and $\o$ is a
positive word in $a_1$ and $a_2$ and the word contains at least
three syllabuses (and ends with $a_2$).
 Then $l(\d_4)<0$.

Suppose that the claim is not true.
Then we have   $l(\o)=0$, $l(\d_1)=0$ and $\d_1$
is a word as in Case (B1) or (B2) or (B3).
But in each of such cases for $\d_1$, one can easily check that
$l(\o)=0$ cannot not hold.
This contradiction proves Claim 6.

We are now ready to finish the proof of  Lemma \ref{lk1}. We first
consider the case when $q=0$, i.e. $\b=\d\in P^a$ is an index-3
reduced positive word containing at least four syllabuses and
$n_\sd>0$. If $n_\sd=1$, then the only syllabus in $a_3$ occurring
in $\d$ must belong to a portion of $\d$ that looks like either in
Claim 3 or in Claim 5, since $\d$ contains at least four
syllabuses. Hence,  Lemma \ref{lk1}  holds in this case by Claim 3
and Claim 5. If $n_\sd\geq 2$, then
 Lemma \ref{lk1} follows from  Claims 3 and 4.

We now consider the case when $q=1$, i.e.
$\b=a_i^{-1}\d\in P^a$ is cyclically reduced,
where
$\d \in P^*$ is index-3 reduced with $n_\sd>0$ and contains
at least four syllabuses.
Suppose otherwise that $lk(L)>0$.
We will get a contradiction.
We have three subcases to consider, corresponding to $i=1,2,3$.
Note that by  Claims 3-5, we have
$l(\d)\leq 0$.

If $i=1$ (i.e. $\b=a_1^{-1}\d$),
then we must have
$l(\d)=0$ and $l(a_1^{-1})=1$,
and the first strand of $\b$ is from $L_1$ and the second
strand of $\b$ from $L_2$.
Also,  $\d$ must start with a syllabus in $a_2$ or $a_3$,
and end with a syllabus in $a_2$ or $a_3$.

Suppose that $\d$ starts with syllabus $a_2^k$ and ends with
syllabus $a_3^j$. Then $\b=a_1^{-1}a_2^k\cdots  a_1^pa_2^ma_3^j$,
which  is conjugate to $\b'=a_1^{-1}a_2^{k+j}\cdots  a_1^pa_2^m$.
The braid $\b'$ is still in $P^a$ but has one less number of
syllabuses in $a_3$.
 So if $\b'$ contains at least five syllabuses, we may use
 induction on $n_\sd$ to conclude that Lemma \ref{lk1} holds in this
 case.
 We may then assume that
 $\b'$ contains four syllabuses and, thus,
 $\b=a_1^{-1}a_2^ka_1^pa_2^ma_3^j$.
 But by Claim 5, $\d$ is as in case (C3);
the strand 1 of $\d$ is from $L_1$ and the strands 2 and 3 of $\d$
are from $L_2$. But this does not match with
the link component  assignment already given to
the strands of $a_1^{-1}$.

Suppose that $\d$ starts with syllabus $a_2^k$ and ends with
syllabus $a_2^j$. Then each syllabus in $a_3$ is followed by a
word without $a_3$ in which both $a_1$ and $a_2$ occur. So we may
decompose $\b$ into portions as $\b=a_1^{-1}\eta_0\eta_1\cdots
\eta_n$ such that $\eta_0\in P^a$ is a word without $a_3$ and each
of $\eta_1,...,\eta_n\in P^a$ is a word that starts with a
syllabus in $a_3$ which is followed by a word without $a_3$ but
containing both $a_1$ and $a_2$. It follows from Claim 3 we must
have $l(\eta_0)=l(\eta_1)=\cdots  =l(\eta_n)=0$ and each of
$\eta_1$,..., $\eta_n$ is as in one of the situations described in
(B1)-(B3). But any two cases of (B1)-(B3) will not match with
their link component assignments if they are adjacent portions of
$\d$. Hence $n=1$. By Claim 6, we have $\eta_0=a_2^k$. One can now
easily check for $\b=a_1^{-1}a_2^k\eta_1$ that the conditions
$l(a_1^{-1})=1$, $l(a_2^k)=0$ and $l(\eta_1)=0$ (so $\eta_1$ is as
in one of (B1)-(B3)) imply that there is no consistent link
component assignment to the strands of $a_1^{-1}$, $a_2^k$ and
$\eta_1$.

Similarly as in the previous case, we can get a contradiction for
the case when
 $\d$ starts with syllabus $a_3^k$ and ends with
syllabus $a_2^j$.

Suppose that $\d$ starts with syllabus $a_3^k$ and ends with
syllabus $a_3^j$. Then $\b=a_1^{-1}a_3^k\cdots  a_1^pa_2^ma_3^j$,
which is conjugate to $\b'=a_1^{-1}a_2^ja_3^k\cdots a_1^pa_2^m$.
The braid $\b'$ is still in $P^a$ but has one less number of
syllabuses in $a_3$.
 Also $\b'$ contains at least five syllabuses, so we may use
 induction to see that Lemma \ref{lk1} holds in this
 case.

Similarly one can deal with the case that $i=2$.

Finally, we consider the case that $i=3$, i.e. $\b=a_3^{-1}\d$.
Note that $\d$ cannot start or end with $a_3$.

Suppose that  $\d$ starts with $a_1$ and ends with $a_2$. Then
$\b=a_3^{-1}a_1^ja_2^k\cdots a_1^ma_2^n$ for some $j,k,m,n>0$.
Using the identity $a_3^{-1}a_1=a_2a_3^{-1}$, we get
$\b=a_2^ja_3^{-1}a_2^k\cdots a_1^ma_2^n
=a_2^ja_2a_1^{-1}a_2^{k-1}\cdots a_1^ma_2^n$. So $\b$ is conjugate
to $\b'=a_1^{-1}a_2^{k-1}\cdots a_1^ma_2^{n+j+1}$. Obviously
$\b'\in P^a$ unless $\b=a_3^{-1}a_1^ja_2a_1^pa_2^q\cdots
a_1^ma_2^n$. If $\b=a_3^{-1}a_1^ja_2a_1^pa_2^q\cdots a_1^ma_2^n$,
then it conjugates to $\b''=a_1^pa_2^q\cdots a_1^ma_2^n$ and thus
$l(\b)\leq 0$. Hence we may assume that $\b\in P^a$. So if $\b'$
contains at least five syllabuses, then we are back to the treated
case $i=1$. Note that if $k>1$, then $s_{\sb'}>4$. So we may
assume that $k=1$ and that $\b'=a_1^{-1}a_3^pa_1^ma_2^{n+j+1}$.
Applying  Claim 3, we see that $l(\b)=l(\b')\leq 0$.

Suppose that  $\d$ starts with $a_2$ and ends with $a_2$. Then
$\b=a_3^{-1}a_2^k\cdots a_1^ma_2^n$$=a_2a_1^{-1}a_2^{k-1}\cdots
a_1^ma_2^n$ for some $k,m,n>0$. So $\b$ is conjugate to
$\b'=a_1^{-1}a_2^{k-1}\cdots a_1^ma_2^{n+1}$. Again using Claim 3,
we see that  $l(\b)=l(\b')\leq 0$.

Similarly one can deal with the case that
  $\d$ starts with $a_1$ and ends with $a_1$, and the case
  that $\d$ starts with $a_2$ and ends with $a_1$.
 The proof of Lemma
\ref{lk1} is now finished.
 \qed

We are now ready to prove Proposition \ref{casson}, which is the
content of the rest of this section. Let $\b=a_i^{-q}\d$ be a
3-braid in $P^a$  as given in the proposition,  whose closure is
the given knot $K$. If some  syllabus in $a_3$ occurring in $\d$
has power $k>2$, we apply the crossing change formula $(*)$ for
Casson invariant to $K=\hat\b$ at the second crossing of the
syllabus. The link $L$ (of two components) obtained by smoothing
the crossing is the closure of a braid of the type as described in
Lemma \ref{lk1} and thus  has non-positive linking number. So  we
get a new braid $\b'$ which is identical with $\b$ except with two
less in power at the syllabus in $a_3$ and the  Casson invariant
of $\hat \b'$ is less than or equal to that of $\hat\b$. Obviously
$\b'$ is still in $P^a$. So it suffices to show that $\hat\b'$ has
positive Casson invariant if it is not in the excluded set $E$
given in Proposition \ref{casson}. If $\b'$ is in the set $E$,
then one can verify directly that the original braid $\b$ has
positive Casson invariant. Similarly we may reduce the power of a
syllabus in $a_2$ or in $a_1$ to one or two,  using the formula
$(*)$,  so that the resulting new braid is still in $P^a$, without
increasing the Casson invariant, and that if the new braid is in
the set $E$, then the old braid has positive Casson invariant.
Therefore we only need to show the proposition under the extra
condition that every syllabus of $\d$ has power one or two. We
shall prove this by induction on the number $n_\sd$.

We first consider the initial step of the induction, i.e. the case
 when $n_\sd=0$.
 If $q=0$,
then $K=\hat \b=\hat\d$ is a positive knot. So by  the proof of
Theorem \ref{thm1},  Proposition \ref{casson} holds in this case.
Actually in this case one can easily give a quick self-contained
proof as follows.
 First apply
the formula $(*)$ to any syllabus of  $\b$ which has power larger
than one so that  $\b$ becomes a new braid  with two less
crossings and the Casson invariant of $K$ is equal to the Casson
invariant of the closure of the new braid minus some negative
 integer (by Lemma \ref{lk0}).
So after several such steps, the braid  $\b$ is simplified to a
braid $\b_0$ in which every syllabus has power one such that the
Casson invariant of $K=\hat\b$ is larger than that of $\hat \b_0$
unless $\b=\b_0$. But $\hat \b_0$ is a $(3,n)$ torus knot. The
normalized Alexander polynomial of a $(3,n)$-torus knot $T(3,n)$
is $$\D_{T(3,n)}(t)=\frac{(t^{3n}-1)(t-1)}{(t^3-1)(t^n-1)t^{n-1}}.
$$ Pure calculation of the second derivative of $\D_{T(3,n)}(t)$
valued at $t=1$ gives $$\frac12\D_{T(3,n)}''(1)=\frac{n^2-1}{3}.$$
The conclusion of Proposition \ref{casson} follows in this case.

Suppose then that $q=1$ (and $n_\sd=0$). If $\b$ contains exactly
five or six syllabuses and at least one of them has power $2$ then
one can verify directly that $\hat\b$ has positive Casson
invariant (the checking is pretty quick
 since every syllabus of $\b$ has power
 at most two and also note that $\b$ is
cyclically reduced and $\hat \b$ is a knot). If $\b$ has more than
six syllabuses and one of them has power $2$, then we apply the
formula $(*)$ at such syllabus to reduce the number of syllabuses
of $\b$. If the new braid $\b'$ is not in $P^a$, then a
cancellation must occur between $a_i^{-1}$ and an $a_i$ in $\b'$.
After the cancellation, we get a braid which is a positive word in
$a_1$ and $a_2$ and thus its closure has positive Casson invariant
unless it is a trivial knot. So Lemma \ref{lk} implies that the
Casson invariant of the old knot $\hat \b$ is positive. If the new
braid $\b'$
 contains at least seven syllabuses,  then we may continue to do
such simplification. Note that $s_{\sb'}$ is one or two less than
$s_\sb$ when $\b'$ is still in $P^a$. So we may assume that every
syllabus of $\d$ has power one. Hence $\d$ looks like
$\d=(a_1a_2)^p$ or $\d=(a_1a_2)^pa_1$ or $\d=(a_2a_1)^p$ or
$\d=(a_2a_1)^pa_2$ for some $p\geq 2$. The case $\d=(a_1a_2)^p$ or
the case $\d=(a_2a_1)^p$ cannot happen since in such case we must
have  $i=3$ and the closure of $\b$ is then not a knot. Consider
the case $\d=(a_1a_2)^pa_1$. We have  $i=2$ or $3$. If $i=2$, then
to be a knot, $\b=a_2^{-1}(a_1a_2)^pa_1$ for $p>2$ and $p\ne 2$
(mod $3$). Now one can  verify directly that $C_{\hat\sb}>0$. If
$i=3$,  then $\b=a_3^{-1}(a_1a_2)^pa_1$ is conjugate to
$a_1(a_1a_2)^{p-1}a_1$ which obviously has positive Casson
invariant. The case  $\d=(a_2a_1)^pa_2$ can be treated similarly.
The proof of  Proposition \ref{casson} for the initial step
$n_\sd=0$ is complete.

Now we may assume that $n_\sd>0$. We warn the reader that the rest
of the proof will involve more patient case counting. Nevertheless
 the guiding  idea will still be more or less as follows: apply
the crossing change formula $(*)$ for the Casson invariant and
Lemma \ref{lk1} at a suitable chosen crossing of  the given braid
$\b=a_i^{-q}\d$ to reduce its complexity, simplify the new braid
(if necessary) to get a braid $\b'=a_i^{-q}\d'$ in $P^a$, use the
induction if $s_{\sd'}\geq 4$ and $\b'$ is not in the excluded set
$$E=\left\{\begin{array}{llll} a_1^{-1}a_2a_3^2a_1a_2,
&a_1^{-1}a_3^2a_1a_2a_3,
&a_1^{-1}a_3a_1a_2^2a_3,&a_1^{-1}a_2a_3a_1a_2^2,\\
a_2^{-1}a_3a_1a_2a_3^2,&a_2^{-1}a_3a_1^2a_2a_3,
&a_2^{-1}a_1a_2a_3^2a_1, &a_2^{-1}a_1^2a_2a_3a_1,\\
a_3^{-1}a_1a_2a_3a_1^2,&a_3^{-1}a_1a_2^2a_3a_1,
&a_3^{-1}a_2a_3a_1^2a_2, &a_3^{-1}a_2^2a_3a_1a_2
 \end{array} \right
\},$$ otherwise verify directly that the original braid $\b$ has
positive Casson invariant. Besides we shall often make use of the
 conditions such as $\hat\b$ is a knot, $\b$ is cyclically
reduced and $\d$ is index-3 reduced.

\noindent {\bf Claim  D1}. If
  $a_3^2$ is a syllabus of $\d$,  then the Proposition holds.

Consider the first such syllabus occurring in $\d$. Applying the
formula $(*)$ to $\hat\b$ at the second crossing of the syllabus
$a_3^2$, we get a new braid $\b'=a_i^{-q}\d'$, where $\d'$ is the
braid obtained from $\d$ by deleting the syllabus $a_3^2$,  such
that $\d'\in P^*$, $n_{\sd'}=n_\sd-1$, and $C_{\hat\sb}\geq
C_{\hat\sb'}$ (by Lemma \ref{lk1}). If $\b'$ is still in $P^a$ and
$\d'$ contains at least four syllabuses but $\b'$ is not a word in
the excluded set $E$, then we may use induction.

Suppose that $\b'\in E$. Since the syllabus we are considering is
the  first such occurring in $\b$ and since $\d$ is index-3
reduced,  $\b$  can only be the word $a_2^{-1}a_3^2a_1a_2a_3^2a_1$
or $a_2^{-1}a_3^2a_1^2a_2a_3a_1$  or $a_2^{-1}a_3^2a_1a_2a_3a_1^2$
or $a_1^{-1}a_2a_3a_1a_2^2a_3^2$. In such a case one can verify
directly that $C_{\hat \sb}$ is positive.

Hence, we may assume that either $\b'$ is  still in $P^a$ with
$\d'$ containing exactly three syllabuses, or $\b'$ is no longer
in $P^a$. In the former case,  $\b$ is a word in the set
$$\left\{\begin{array} {lllll}a_1^ja_2^ka_3^2a_1^m,&
a_2^ja_1^ka_2^ma_3^2,& a_2^ja_3^2a_1^ka_2^m,&
a_3^2a_1^ka_2^ja_3^m,& a_3^2a_1^ka_2^ja_1^m,\\
a_3a_1^ka_2^ma_3^2,& a_1^{-1}a_2^ja_1^ka_2^ma_3^2,&
a_1^{-1}a_2^ja_3^2a_1^ka_2^m,& a_1^{-1}a_3a_1^ka_2^ma_3^2,&
a_2^{-1}a_1^ja_2^ka_3^2a_1^m,\\a_2^{-1}a_3^2a_1^ja_2^ka_3^m,&
a_2^{-1}a_3^2a_1^ka_2^ja_1^m,& a_3^{-1}a_1^ja_2^ka_3^2a_1^m,&
a_3^{-1}a_2^ja_3^2a_1^ka_2^m\end{array}\right\}$$ for some
$j,k,m\in\{1,2\}$. When $\b$ is one of words in this set but is
not in the excluded set $E$, one can
 verify directly using
the formula $(*)$ that $\hat\b$ has positive Casson invariant.
 For instance, when
 $\b=a_1^{-1}a_2^ja_1^ka_2^ma_3^2$, it is conjugate to
$a_2^ja_1^ka_2^ma_3^2a_1^{-1}=a_2^ja_1^ka_2^ma_1^{-1}a_2^2$ which
in turn is conjugate to $a_1^{-1}a_2^{j+2}a_1^ka_2^m$. So we need
to show that the Casson invariant of the closure of the braid
$\eta=a_1^{-1}a_2^{j+2}a_1^ka_2^m$ is positive. We have  eight
possible cases for $\eta$ corresponding to various possible values
of $j$, $k$ and $m$. But only  in case $j=k=m=1$ or case $j=1,
k=m=2$ or case $j=k=2,m=1$,   the closure of the involved braid is
a knot, and in such a case one can verify directly using formula
$(*)$ that $\hat \b=\hat\eta$ has positive Casson invariant. In a
similar way one can verify the proposition for each of the other
words in the above set. (We did the checking!)

We now consider the latter case when  $\b'$ is not in $P^a$. It
follows that the syllabus $a_3^2$ is either the first or the last
syllabus of $\d$ and $q=1$. We consider the case when the syllabus
$a_3^2$ is the first syllabus of $\d$. The case when
 the syllabus $a_3^2$  is the last  syllabus of $\d$
can be treated similarly. It follows that
$\b=a_1^{-1}a_3^2a_1^ja_2^k\cdots  $, for some $j,k\in \{1,2\}$,
and $\b$ does not end with $a_1$. If $j=2$,  then
$\b'=a_1^{-1}a_1^2a_2^k\cdots $ which is isotopic to $\b''
=a_1a_2^k\cdots  $ which is in $P^a$. So if $\b''$ contains at
least four syllabuses, we may use the induction.  We may then
assume that $\b''$ has less than four syllabuses. It follows that
$\b''=a_1a_2^ka_3^m$ which is a knot only when $k=2, m=1$ or
$k=1,m=2$  and in these two cases $C_{\hat\sb}\geq
C_{\hat\sb''}>0$. So we may assume that $j=1$. In this case $\b'$
is isotopic to $\b''=a_2^k\cdots  $ which  is in $P^a$.  Hence if
$\b''$ contains at least four syllabuses, we may use the
induction.  If $\b''$ has less than four syllabuses,   then we
must either have
 $\b= a_1^{-1}a_3^2a_1a_2^ka_1^na_2^m$ or $\b=
a_1^{-1}a_3^2a_1a_2^ma_3^j$. In the former case, $\b''$ has
positive Casson invariant (a nontrivial positive knot). In the
latter case, only when $m=1, j=1$, $\hat\b$ is a knot. But this
braid is in the  set $E$.
 The proof of Claim D1 is now complete.

By Claim D1, we may now assume that every syllabus in $a_3$
occurring in $\d$ has power equal to one.

\noindent {\bf Claim D2}. We may assume that every syllabus in $a_1$ occurring in
 $\d$ has power equal to one.

Suppose that $\d$ contains syllabuses in $a_1$ of power two.
 Consider the first such syllabus occurring in $\d$.
Applying the formula $(*)$ to $\hat\b$ at the first crossing of
the syllabus $a_1^2$, we get a new braid $\b'=a_i^{-q}\d'$ such that the
Casson invariant of $\hat\b'$ is less than or equal to that of $\hat\b$ by Lemma
\ref{lk1}. Obviously $\d'$ is still a positive word in
$a_1,a_2,a_3$ but may not be in $P^*$ or $P^a$. We have several
possibilities for $\d$ around the given syllabus $a_1^2$:
$\d=\cdots  a_2^ja_1^2a_2^k\cdots  $ or $\d=\cdots  a_3a_1^2a_2^k\cdots  $ or
$\d=a_1^2a_2^k\cdots  $ or $\d=\cdots  a_1^2$, for some $j,k\in \{1,2\}$.

\noindent Case (D2.1).  $\d=\cdots  a_2^ja_1^2a_2^k\cdots  $.

Then $\d'=\cdots  a_2^{j+k}\cdots  $ and $\b'=a_i^{-q}\d'$ is still in
$P^a$. Also if $j+k>2$, we may apply the
formula $(*)$ one more time to bring  it down to one or two.
Let $\b''=a_i^{-q}\d''$ be the resulting braid.
Then if $\d''$ contains at least four syllabuses and $\b''$ is not in $E$,
then we
have eliminated one $a_1^2$.
If $\b''$ is in $E$,  then one can verify directly that
the original braid $\b$ has positive Casson invariant.
 Note
that $s_{\sd'}=s_{\sd}-2$. Suppose that $s_{\sd'}<4$. Then $s_\sd$
is four or five and $\b=a_i^{-q}a_2^ja_1^2a_2^ka_3$ or
$\b=a_i^{-q}a_1a_2^ja_1^2a_2^ka_3$ or
$\b=a_i^{-q}a_3a_1a_2^ja_1^2a_2^k$. Easy to check that when $q=0$,
any knot from these cases has  positive Casson invariant.
So assume that $q=1$ in these cases.

If $\b=a_i^{-1}a_2^ja_1^2a_2^ka_3$,  then $i=1$ since $\b$ is
cyclically reduced. So $\b=a_1^{-1}a_2^ja_1^2a_2^ka_3$. To be a
knot, we have $k=j=1$ or $k=j=2$. In each of the two cases, one
can verify directly that the Casson invariant of the knot is
positive.

If $\b=a_i^{-q}a_1a_2^ja_1^2a_2^ka_3$,  then $i=2$. That is
$\b=a_2^{-1}a_1a_2^ja_1^2a_2^ka_3$. One can also directly verify
that $C_{\hat\sb}>0$ (for those values of $k,j\in\{1, 2\}$
 which make $\hat \b$ a knot).

The case that  $\b=a_i^{-q}a_3a_1a_2^ja_1^2a_2^k$ can be treated
similarly.

\noindent Case (D2.2).  $\d=\cdots  a_3a_1^2a_2^k\cdots  $.

Then $\d'=\cdots  a_3a_2^k\cdots  =\cdots  a_2a_1a_2^{k-1}\cdots  $.

\noindent Case (D2.2.1.)  $k=2$.

If $\d$ does not start with $a_3$,  then $\b'=a_i^{-q}\d' =
a_i^{-q}\cdots  a_2^{j+1}a_1a_2\cdots  $ is still in $P^a$ and is not in the set
$E$. In such case if
$\d'$ contains at least four syllabuses, we may apply induction
since $n_{\sd'}=n_\sd-1$.
If $s_{\sd'}<4$, then
$\b=a_i^{-q}a_2^ja_3a_1^2a_2^2$. To be a knot, we have  $q=0, j=1$
or $q=1,i=1,j=2$. In each of the two cases  we have  $C_{\hat\sb}>0$ by
direct calculation.

If  $\d$ starts with $a_3$ but $q=0$,  then $\b'=\d'$ is still in $P^a$
but not in $E$. Also
$s_\sb=s_\sd=s_{\sb'}=s_{\sd'}$ and $n_{\sd'}=s_\sd-1$.
So we may apply induction.

If  $\d$ starts with $a_3$ and $q=1$,  then $i=1$ or $2$. If $i=1$,
then $\b'=a_1^{-1}a_2a_1a_2\cdots  $ is still in $P^a$ but not in $E$
and contains at least
five syllabuses. So we may use
induction since $s_{\sd'}=s_\sd$ but $n_{\sd'}=n\sd-1$. If $i=2$,
then $\b'=a_1a_2\cdots  $ is in $P^a$ and is not in $E$.
Also $s_{\sb'}
=s_\sd-1$. Hence if $s_\sd>4$, we may use
induction. So
suppose that $s_\sd=4$. Then $\b=a_2^{-1}a_3a_1^2a_2^2a_3$ or
$\b=a_2^{-1}a_3a_1^2a_2^2a_1^j$. But the former is not a knot. The
latter is a knot when $j=2$, in which case $C_{\hat\sb}>0$.

\noindent Case (D2.2.2). $k=1$.

Then  $\d=\cdots  a_3a_1^2a_2a_3\cdots  $ or $\d=\cdots  a_3a_1^2a_2a_1^j\cdots  $ or
$\d=\cdots  a_3a_1^2a_2$. Correspondingly, we have
 $\d'=\cdots  a_2a_1a_3\cdots  =\cdots  a_2^2a_1\cdots  $
  or $\d=\cdots  a_2a_1^{j+1}\cdots  $ or
$\d=\cdots  a_2a_1$. If  $\b'$ is in $P^a-E$ and $s_{\sd'}\geq 4$,
we may use induction since $n_{\sd'}<n_{\sd}$. If $\b'\in E$, then
one can check that the old braid $\b$ has positive Casson
invariant. If $\b'\in P^a$ but $s_{\sd'}<4$, then one can also
verify that  $\b$ always has positive Casson invariant, applying
the conditions that (1) $\b\in P^a-E$, (2) $s_\sd\geq 4$ and (3)
$\hat\b$ is a knot.

\noindent Case (D2.3).  $\d=a_1^2a_2^k\cdots  $.

Then $\d=a_1^2a_2^ka_1^ja_2^m\cdots  $ or
$\d=a_1^2a_2^ka_3a_1^j\cdots  $. And $\d'=a_2^ka_1^ja_2^m\cdots  $
or $\d'=a_2^ka_3a_1^j\cdots  $. So $\b'$ is in $P^a$ unless $\b$
starts with $a_2^{-1}$. If $\b$ starts with $a_2^{-1}$,  then
$\b'= a_2^{k-1}a_1^ja_2^m\cdots  $ or $\b'=a_2^{k-1}a_3a_1^j\cdots
$ which is in $P^a$. Again in each of these cases, if $\b'$ is the
set $E$ or $s_{\sd'}<4$, one can calculate directly that the old
braid $\b$ has positive Casson invariant. Otherwise one can use
the  induction.

\noindent Case (D2.4).  $\d=\cdots  a_1^2$.

This case can be treated similarly as in the previous case.

So by Claims D1 and D2, we  now assume that every  syllabus in
$a_3$  and in $a_1$ occurring in $\d$
 have power equal to one.

\noindent {\bf Claim D3}.  We may assume that every syllabus in $a_2$ occurring in
 $\d$ has power equal to one.

This  can be proved similarly as Claim D2.

So we  now assume that every syllabus  occurring in $\d$
 has power one.

 If $a_1a_2a_1$ appears  immediately after an $a_3$,  then
$\b=\cdots  a_3a_1a_2a_1\cdots   =\cdots  a_3a_1a_3a_2\cdots   =\cdots
a_3a_2a_1a_2\cdots  =\cdots  a_2a_1^2a_2\cdots  $ and so we get an
isotopic braid
in $P^a$ with less number of syllabuses in $a_3$. Hence we may
apply induction unless $\b$ is in the set
\{$a_3a_1a_2a_1$, $a_2^{-1}a_3a_1a_2a_1$, $a_3^{-1}a_2a_3a_1a_2a_1$,
$a_1^{-1}a_3a_1a_2a_1a_2$, $a_1^{-1}a_2a_3a_1a_2a_1a_2$,
$a_3^{-1}a_2a_3a_1a_2a_1a_2$\}.
If $\b$ is a word in this set, then either $\hat\b$ is not a knot or
$\hat \b$ has positive Casson invariant.
 So we may assume that no  $a_3a_1a_2a_1$
occurs in $\b$. A similar argument shows that we may assume that
no $a_2a_1a_2a_3$  occurs in $\b$.
Hence
we may assume that $\b$ is one of the words in the set
{\small $$\left\{\begin{array}
{llll}
a_1a_2a_3a_1,&
(a_3a_1a_2)^m,& (a_3a_1a_2)^ma_3, &(a_3a_1a_2)^ma_3a_1, \\
a_2(a_3a_1a_2)^m,&
a_1a_2(a_3a_1a_2)^m,&
a_2(a_3a_1a_2)^ma_3, &
a_2(a_3a_1a_2)^ma_3a_1,\\
a_1a_2(a_3a_1a_2)^ma_3, &a_1a_2(a_3a_1a_2)^ma_3a_1,&
a_2^{-1}a_1a_2a_3a_1,& a_3^{-1}a_1a_2a_3a_1,\\
a_1^{-1}(a_3a_1a_2)^m,&
a_1^{-1}(a_3a_1a_2)^ma_3,&a_2^{-1}(a_3a_1a_2)^ma_3,&
a_2^{-1}(a_3a_1a_2)^ma_3a_1\\

 a_1^{-1}a_2(a_3a_1a_2)^m,&
a_3^{-1}a_2(a_3a_1a_2)^m,& a_3^{-1}a_1a_2(a_3a_1a_2)^m,&
a_1^{-1}a_2(a_3a_1a_2)^ma_3,\\a_3^{-1}a_2(a_3a_1a_2)^ma_3a_1
&a_2^{-1}a_1a_2(a_3a_1a_2)^ma_3, &

a_2^{-1}a_1a_2(a_3a_1a_2)^ma_3a_1,&
a_3^{-1}a_1a_2(a_3a_1a_2)^ma_3a_1

\end{array}\right\}$$}
where  $m>0$.

If $\b=a_1a_2a_3a_1$, then it has positive Casson invariant.
The case $\b=(a_3a_1a_2)^m$ cannot occur since its closure is not
a knot for all $m>0$.
Similarly each of the cases  $a_2(a_3a_1a_2)^ma_3a_1$, $a_1a_2(a_3a_1a_2)^ma_3$,
$a_2^{-1}a_1a_2a_3a_1$,
$a_3^{-1}a_1a_2a_3a_1$,  $a_1^{-1}(a_3a_1a_2)^m$, $a_2^{-1}(a_3a_1a_2)^ma_3a_1$
and $a_1^{-1}a_2(a_3a_1a_2)^ma_3$ cannot occur as $\b$.
If $\b=(a_3a_1a_2)^ma_3$, then it is conjugate to
$\b'=a^2_3a_1a_2(a_3a_1a_2)^{m-1}$.
So we may apply Claim D1 (note that $n_{\sb'}<n_\sb$) unless $m=1$.
But when $m=1$, one can calculate directly that
 $\b=a_3a_1a_2a_3$ has positive Casson invariant.
 If $\b=(a_3a_1a_2)^ma_3a_1$, then it is conjugate to
$\b'=a_2a_1^3a_2(a_3a_1a_2)^{m-1}$. So we are back to a previous
case (note that $n_{\sb'}<n_\sb$) unless $m=1$. But when $m=1$,
$\hat \b$ is not a knot. If $\b=a_2(a_3a_1a_2)^m$, then it is
conjugate to $\b'=(a_3a_1a_2)^{m-1}a_3a_1a_2^2$. So we may apply
Claim D2 (note that $n_{\sb'}=n_\sb$) unless $m=1$. But when
$m=1$, one can calculate directly that
 $\b=a_2a_3a_1a_2$ has positive Casson invariant.
 Similarly one can deal with the cases
 $a_1a_2(a_3a_1a_2)^m$, $a_2(a_3a_1a_2)^ma_3$
 and  $a_1a_2(a_3a_1a_2)^ma_3a_1$.

Suppose that  $\b=a_1^{-1}(a_3a_1a_2)^ma_3$. To be a knot,
we have $m\geq 3$ and $m=3 (mod 2)$. Also $\b$ is conjugate to
$\b'=a_1^{-1}a_2(a_3a_1a_2)^m$ which has less number of
syllabuses in $a_3$ and thus we may apply the induction.
Similarly one deal with the case when $\b=a_2^{-1}(a_3a_1a_2)^ma_3$.

If $\b=a_1^{-1}a_2(a_3a_1a_2)^m$,   then to be a knot $m$ must be
even. Applying the formula $(*)$ at the first crossing of $\b$, we
get $C_{\hat \sb}= C_{\hat\sb_1}+lk(\hat \l_1)$ where
$\b_1=a_1a_2(a_3a_1a_2)^m$ and $\l_1=a_2(a_3a_1a_2)^m$. One can
easily deduce from the projection of $\l_1$ that $lk(\hat
\l_1)=-m/2$. Hence, we have $C_{\hat\sb}=C_{\hat\sb_1}-m/2.$ In
such case it suffices  to show the following.

\noindent {\bf Claim D4}. $C_{\hat\sb_1}>m/2$.

We knew that  $m=2p$ with $p>0$. We shall prove the claim by
induction on the number $p$. Write $\b_1$ as $\b_1=
a_1a_2a_1^{-1}a_2a_1^2a_2a_1^{-1}a_2a_1^2a_2(a_3a_1a_2)^{m-2}$.
Applying formula $(*)$ to $\b_1$ at the first crossing of the
first syllabus $a_1^2$, we get $C_{\hat\sb_1}=
C_{\hat\sb_2}-lk(\hat\l_2)$ where
$\b_2=a_1a_2a_1^{-1}a_2^2a_1^{-1}a_2a_1^2a_2(a_3a_1a_2)^{m-2}$ and
$\l_2=a_1a_2a_1^{-1}a_2a_1a_2a_1^{-1}a_2a_1^2a_2(a_3a_1a_2)^{m-2}$.
One can easily deduce from the projection of $\l_2$ that $lk(\hat
\l_2)=-2(m-2)-3=-4p+1$. Hence, we have
$C_{\hat\sb_1}=C_{\hat\sb_2}+4p-1.$ We then apply formula $(*)$ to
$\b_2$ at the first crossing of the first syllabus $a_2^2$, we get
$C_{\hat\sb_2}= C_{\hat\sb_3}-lk(\hat\l_3)$ where
$\b_3=a_1a_2a_1^{-2}a_2a_1^2a_2(a_3a_1a_2)^{m-2}$ and
$\l_3=a_1a_2a_1^{-1}a_2a_1^{-1}a_2a_1^2a_2(a_3a_1a_2)^{m-2}$. One
can calculate to see  that $lk(\hat \l_3)=-p-1$. Hence, we have
$C_{\hat\sb_2}=C_{\hat\sb_3}+p+1.$ We then apply formula $(*)$ to
$\b_3$ at the first crossing of the syllabus $a_1^{-2}$, we get
$C_{\hat\sb_3}= C_{\hat\sb_4}+lk(\hat\l_4)$ where
$\b_4=a_1a_2^2a_1^2a_2(a_3a_1a_2)^{m-2}$ and
$\l_4=a_1a_2a_1^{-1}a_2a_1^2a_2(a_3a_1a_2)^{m-2}$. Also one can
calculate to see that $lk(\hat \l_4)=-p$. Hence, we have
$C_{\hat\sb_3}=C_{\hat\sb_4}-p.$ We then apply formula $(*)$ to
$\b_4$ at the first crossing of the syllabus $a_2^{2}$, we get
$C_{\hat\sb_4}= C_{\hat\sb_5}-lk(\hat\l_5)$ where
$\b_5=a_1^3a_2(a_3a_1a_2)^{m-2}$ and
$\l_5=a_1a_2a_1^2a_2(a_3a_1a_2)^{m-2}$. We also have $lk(\hat
\l_5)=-4(p-1)-2$. Hence we have
$C_{\hat\sb_4}=C_{\hat\sb_5}+4(p-1)+2.$ Applying formula $(*)$ to
$\b_5$ at the first crossing of the syllabus $a_1^{2}$, we get
$C_{\hat\sb_5}= C_{\hat\sb_6}-lk(\hat\l_6)$ where
$\b_6=a_1a_2(a_3a_1a_2)^{m-2}$ and
$\l_6=a_1^2a_2(a_3a_1a_2)^{m-2}$. We also have $lk(\hat
\l_6)=-(p-1)-1=-p$. Hence we have $C_{\hat\sb_5}=C_{\hat\sb_6}+p.$
In summary, we get $$C_{\hat\sb_1}=C_{\hat\sb_6}+9p-2.$$ Now one
can easily see that the claim follows.

Similarly we can treat the rest of cases.
The proof of Proposition \ref{casson} is now complete.

\begin{center}
{\bf 4. Proof of Proposition \ref{lamin}} \end{center}

Given a 3-braid  $\b$ in letters $a_1, a_2, a_3$, whose closure is
a knot, there is a canonical way to construct a Seifert surface
for $\hat \b$ as follows: in the projection plane we have the
braid diagram in its canonical form,  place three rectangular
disks in the space so that disk 1 lies in the projection plane and
is  on the left hand side of the braid, disk 3 also lies in the
projection plane but on the right hand side of the braid, disk 2
lies perpendicularly above the projection plane, each disk having
one side running parallel to the braid from top to the bottom,
then to  each letter $a_1$ ($a_1^{-1}$) occurring in $\b$  use a
half negatively (positively) twisted band connecting disks 1 and
2, to each letter $a_2$ ($a_2^{-1}$) use  a half negatively
(positively) twisted band connecting disks 2 and 3, and to each
letter $a_3$ ($a_3^{-1}$) use a half negatively (positively)
twisted band connecting disks 1 and 3 (behind disk 2). Figure 7
illustrates such construction for
$\b=a_1a_2a_3a_1^{-1}a_3^{-1}a_2^{-1}$.
\begin{figure}[!ht] {\epsfxsize=2in
\centerline{\epsfbox{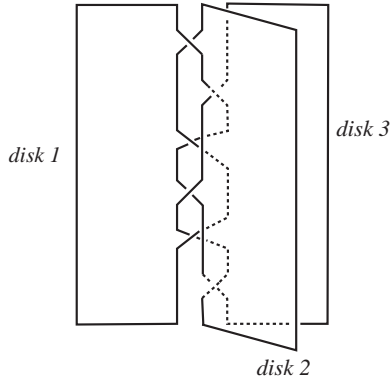}}} \caption{construction of the
canonical Seifert surface }
\end{figure} We call the Seifert
surface of $\hat\b$  so constructed {\it canonical} Seifert
surface of $\hat \b$. In [X], it was proved that if $\b$ is a 3-braid
of norm form (whose definition we recalled in Section 2),  then
its canonical Seifert surface has the  minimal genus (thus is an
essential
and Thurston norm minimizing surface in the exterior of $\hat\b$).

\begin{lemma}\label{lami}
Let $\b$ be a 3-braid such that  $\hat \b$ is a knot. Let $S$ be
the canonical Seifert surface of $\b$ and suppose that it has
minimal genus.
\newline (1) If $S$ contains two half twisted bands corresponding
to the same letter $a_1$ or $a_2$ or $a_3$, then the $-1$-surgery
on $\hat \b$ is a manifold with essential lamination. \newline (2)
If $S$ contains two half twisted bands corresponding to the same
letter $a_1^{-1}$ or $a_2^{-1}$ or $a_3^{-1}$,  then the
$1$-surgery on $\hat \b$ is a manifold with essential lamination.
\end{lemma}

\pf We shall only prove part (1) when $S$ contains two half
twisted bands corresponding to the same letter $a_1$. All other
cases can be proved similarly.

Let  $M$ be the knot  exterior of $\hat \b$ in $S^3$. We shall
also use $M(-1)$ to denote the manifold obtained by Dehn surgery
on the knot $\hat\b$ with the slope $-1$. Let $V$ be the solid
torus filled in $M$ to obtain the manifold $M(-1)$. We first
construct an essential branched surface $B$ in the exterior $M$
and then prove that $B$ (which has boundary on $\p M$) can be
capped off by a branched surface in $V$ to yield  an essential
branched  surface $\hat B$ in $M(-1)$. The construction of $B$ is
similar to that given in [Rs].

\begin{figure}[!ht]
{\epsfxsize=5in \centerline{\epsfbox{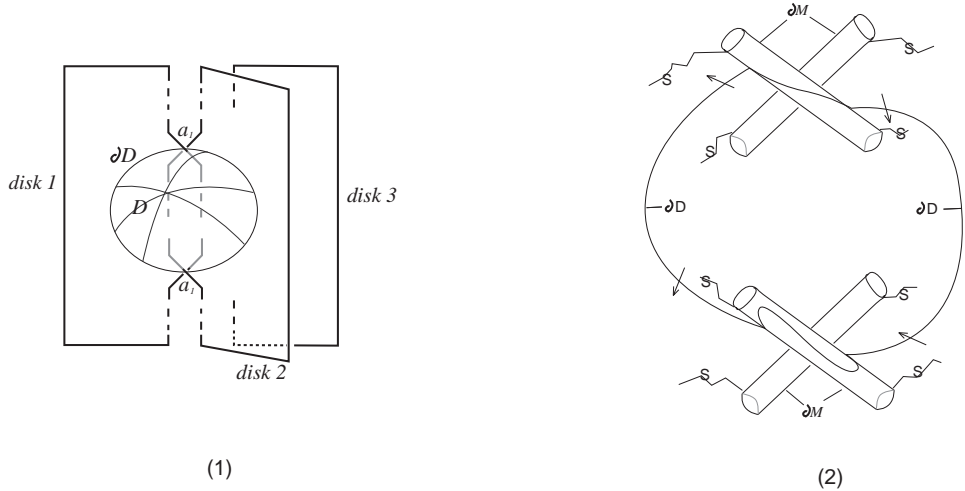}}} \caption{
construction of $B$}
\end{figure}

Since  $S$ contains two half negatively twisted  bands
corresponding to $a_1$,  there is a disk  $D$ in $M$ as shown in
Figure 8 (1) whose boundary lies in $S\cup \p M$ and whose
interior is disjoint from $S\cup \p M$.
 With more detail, the boundary of $\p D$ intersects $S$
in two disjoint arcs and intersects $\p M$ in two disjoint arcs
with the latter happening around the places corresponding to the
two bands of $a_1$. (Similar disks were used in [BMe] for a
different purpose).  The branched surface $B$ is the union of the
Seifert surface $S$ and the disk $D$ with their intersection locus
smoothed as shown in Figure 8 (2). The arrows in the figure
indicate the cusp direction of branched locus.
 An argument as in [Rs] shows that
the branched surface fully carries a lamanition with no compact
leaves and each negative slope can be realized as the boundary
slope of a lamination fully carried by $B$. Since the disk $D$
intersects the knot exactly twice, the branched surface $B$ is
essential in $M$ by  [G2, 3.12].
 The branched locus of $B$ is a set of two disjoint
 arcs properly embedded in $S$, each being non-separating.
The branched surface $B$ meets $\p M$ yielding  a train track in
$\p M$ as show in Figure 9 (1).  Let $\cal L$ be a lamination
fully carried by $B$ whose boundary slope is $-1$. Then $\p \cal
L$ must look like as shown in Figure 9 (2).

\begin{figure}[!ht] {\epsfxsize=5in
\centerline{\epsfbox{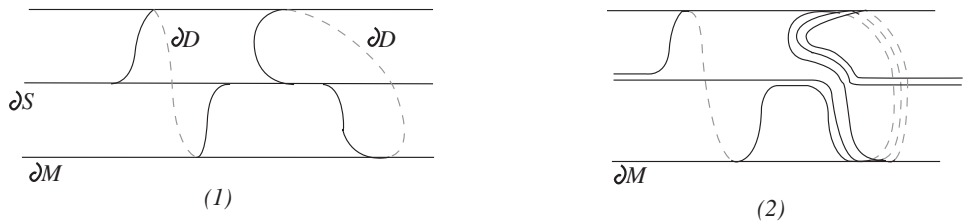}}} \caption{(1) $\p B$ on $\p
M$\hspace{3mm} (2) the curve of slope $-1$ fully carried by $\p
B$}
\end{figure}

We now construct a branched surface $B_{\sv}$ in the sewn solid
torus $V$ such that the train track  $ B_\sv\cap \p V$ is
equivalent to the train track  $B\cap \p M$ and $B_\sv$ fully
carries a lamination which is a set of meridian disks of $V$.
Hence, $B$ and $B_\sv$ match together and form a branched surface
$\hat B$ in $M(-1)$. Take a meridian disk $D_0$ of $V$ and push
part of it near and around $\p V$ as shown in Figure 10 (1) and
then identify two disjoint sub-disks of $D_0$ as shown in Figure
10 (2). This gives a branched surface $B_1$ with the cusp
direction along its  singular locus (an arc) as shown in Figure 10
(2). Then we split $B_1$  locally at a place around a point of $\p
B_1$ as shown in Figure 10 (3) and then we start pinch the
resulting branched surface along part of its boundary
 as shown in Figure 10 (4).  The pinching continues as shown in Figure
 10 (5) until we get  the branch surface whose boundary is as
 shown in Figure 10 (6).
The resulting surface is the branched surface $B_v$. Obviously
 the train track  $\p B_v$ on $\p V$ is equivalent to the train
 track $\p B$ in $\p M$ and they can be matched in $\p V=\p M$.

\begin{figure}[!ht] {\epsfxsize=6in
\centerline{\epsfbox{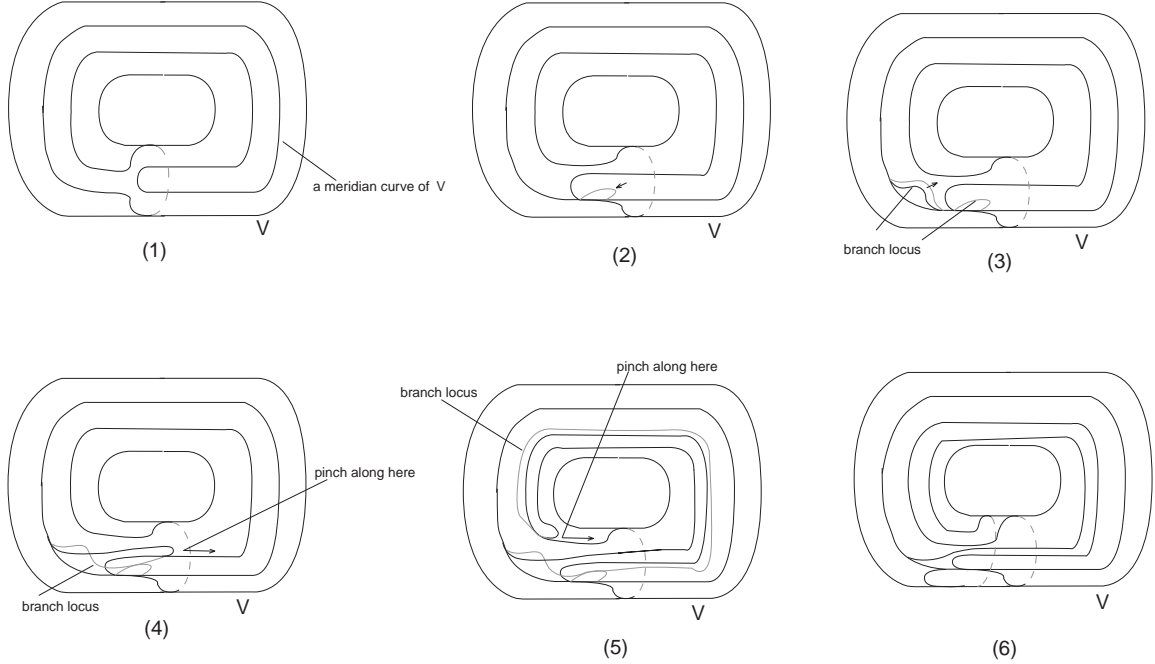}}} \caption{ construction of $B\sv$}
\end{figure}

To see that $\hat B$ is essential
 we have five things to
check  by [GO, Definition 2.1]:\newline (i) $\hat B$ has no discs
of contact;\newline (ii) The horizontal surface $\p_h N(\hat B)$
is incompressible and $\p$-incompressible in
$M(-1)-\stackrel{\circ}{N}(\hat B)$, there are no monogons in
$M(-1)-\stackrel{\circ}{N}(\hat B)$ and no component of $\p_h
N(\hat B)$ is a 2-sphere;
\newline
(iii) $M(-1)-\stackrel{\circ}{N}(\hat B)$ is irreducible;
\newline
(iv) $\hat B$ contains no Reeb branched surface;
\newline
(v) $\hat B$ fully carries a lamination.

Condition (v) follows automatically by the construction since
leaves of a lamination fully carried by $ B$ with boundary slope $-1$ match
on $\p M=\p V$ with (disk) leaves of a lamination
  fully carried by $B_\sv$. It also follows that $\hat B$
does not  carry any compact surface since $B$ does not. Hence in
 particular condition
(iv) holds also for $\hat B$. By the construction, one can easily
see that  $V-\stackrel{\circ}{N}(B_v)$ has two components, each of
which  topologically looks like as shown in Figure 11. It follows
that $M(-1)-\stackrel{\circ}{N}(\hat B)$ is topologically the same
as $M-\stackrel{\circ}{N}(B)$, with the same horizontal surface.
>From this we get conditions (ii) and (iii) for $\hat B$.

We now show that $\hat B$ has no disk of contact. Note that
$\p_v(N(\hat B))$ is a set of two annuli  and each of the annuli
is obtained from matching a component of $\p_v(N(B))$ (a vertical
disk) and a component of $\p_v(N(B_\sv))$ (a vertical disk). Hence,
if $D_c$  were a contact disk in $N(\hat B)$, then its boundary would
have to
intersect a component of
 $\p_v(N(B))$. It follows then that the interior of $D_c$ must enter
 into the region of $N(B)$ correspond to a branch of
 ($S-$ the singular locus of $B$). But one can easily see from Figure 8
 that the complement of the singular locus of $B$ in $S$
 is a connected surface. It follows that $D_c$ has to intersect
 every $I$-fiber of $B$ since $D_c$ is transverse to
 $I$-fibers of $B$. In particular $\p D_c$ has
 to intersect both of the vertical annuli of $\p_v(\hat B_v)$,
 which gives a contradiction.
\qed

\begin{figure}[!ht] {\epsfxsize=3in
\centerline{\epsfbox{3b11.ai}}} \caption{a component of
$V-\stackrel{\circ}{N}(B_v)$}
\end{figure}

We now prove Proposition \ref{lamin}. By [X],
the canonical Seifert surface of $\hat \b$ has
minimal genus.
If the condition (1) of Proposition \ref{lamin} holds,
then  the conclusion of  Proposition \ref{lamin} follows
obviously from Lemma \ref{lami}.
Suppose that  the condition (2) of Proposition \ref{lamin} holds.
To show that the $1$-surgery on $\hat\b$ gives a manifold with
an essential lamination, we may assume, by Lemma \ref{lami},
 that $\eta$ contains at most three syllabus, and they have different
 subscripts and all have  power  $-1$. But $\eta$ contains
at least two  syllabuses. Suppose that the first and the second
syllabuses of $\eta$ are $a_3^{-1}$ and $a_2^{-1}$, i.e.
$\eta=a_3^{-1}a_2^{-1}\cdots  $. Then since $\b$ is a shortest
word, the word $\d$ does not end with a syllabus in $a_3$. Suppose
that $\d$ ends with a syllabus in $a_1$. Then we have $\b=\cdots
a_1a_3^{-1}a_2^{-1}\cdots  $. By a band move isotopy of the
Seifert surface as shown in Figure 12 (1) we get an isotopic
3-braid $\b'$ which contains two $a_2^{-1}$. (Algebraically,
$\b=\cdot a_1a_3^{-1}a_2^{-1}\cdot=\cdot
a_2^{-1}a_1a_2^{-1}\cdot=\b'$).
 Further
the canonical Seifert surface of $\hat \b'$ is isotopic to that of
$\hat\b$ and thus has minimal genus. So we may apply Lemma
\ref{lami} to see that for the knot $\hat \b'=\hat\b$, the
$1$-surgery gives a manifold with essential lamination.
 Suppose then that $\d$ ends with
a syllabus in $a_2$. Since $\d$ is assumed to contain at least two syllabuses,
 $\b=\cdots  a_1a_2^ka_3^{-1}a_2^{-1}\cdots$. Again we
may first use the  band-isotopy as shown in Figure 12 (2) and then
use the band isotopy of Figure 12 (1) to  get an isotopic 3-braid
whose canonical Seifert surface
 contains two bands corresponding to $a_2^{-1}$.
 Hence Proposition \ref{lamin} follows from Lemma \ref{lami}
 in this case as well.
Similarly one can treat the cases when $\eta$ starts with
$a_2^{-1}a_1^{-1}$ or with $a_1^{-1}a_3^{-1}$.
The case when $\d$ contains at most three
syllabuses, each having power at most one, can be proved similarly.
This proves Proposition \ref{lamin} under its condition (2).
Finally if
 the condition (3) of Proposition \ref{lamin} holds then either
 condition (2) of Proposition \ref{lamin} holds or
 one can get directly
 two letters $a_i$ of the same subscript in $\d$ and two letters
 $a_j^{-1}$ of the same subscript  in $\eta$.
 So again the Proposition follows from Lemma \ref{lami}.\qed

\begin{figure}[!ht] {\epsfxsize=5in
\centerline{\epsfbox{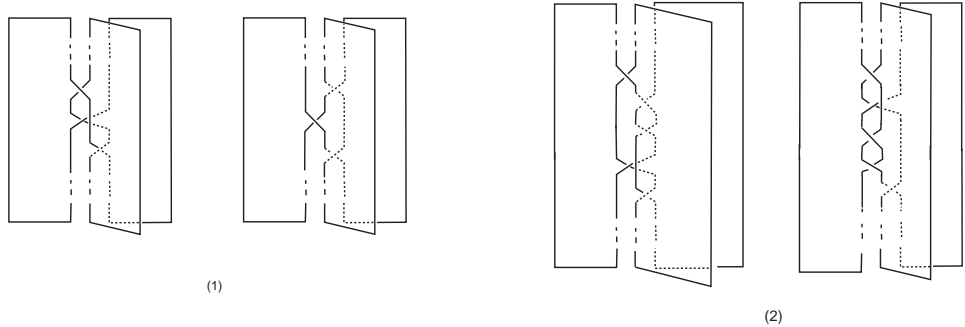}}} \caption{the band move isotopies}
\end{figure}

 \vspace{3mm}
{\small

{\begin{center} {\bf References}
\end{center}

\noindent[AM] S. Akbulut and J. McCarthy, Casson's invariant for
oriented homology 3-spheres. Princeton University Press (1990).

 \noindent
[BMa]  R. Bing and J. Martin,  Cubes with knotted holes.
  Trans. Amer. Math. Soc. 155 (1971)
217--231.

\noindent
[BMe] J. Birman and W. Menasco,
  Studying links via closed braids. III.
  Classifying links which are closed $3$-braids, Pacific
J. Math. 161 (1993) 25--113.

\noindent
[BS]
S. Bleiler and M. Scharlemann,  A projective plane in ${R}\sp 4$ with three critical
points is standard. Strongly invertible knots have property $P$.
Topology 27 (1988) 519--540.

\noindent [CGLS] M. Culler, C. M. Gordon, J. Luecke and P. Shalen,
Dehn surgery on knots. Ann. of Math. 125 (1987) 237-300.

\noindent
[DR] C. Delman and R. Roberts,
Alternating knots satisfy strong property P.
 Comment. Math. Helv. 74 (1999) 376--397.

\noindent [D]  N. Dunfield,  Cyclic surgery, degrees of maps of
character curves, and volume rigidity for hyperbolic manifolds.
 Invent.
Math. 136 (1999) 623--657.

\noindent
[G1]
D.  Gabai, Surgery on knots in solid tori. Topology 28 (1989) 1--6.

\noindent [G2] ---, Foliations and the topology of
$3$-manifolds.  J. Differential Geom. 18 (1983) 445--503.

 \noindent[GO] D. Gabai and U. Oertel, Essential
laminations and 3-manifolds.  Ann. Math. 130 (1989) 41-73.

 \noindent
[GL] C.  Gordon and J. Luecke,  Knots are determined by their complements.
 J. Amer. Math. Soc. 2 (1989), no. 2, 371--415.

 \noindent
 [H] J. Hempel,
 A simply connected $3$-manifold is $S\sp{3}$ if it is the sum of a solid torus
and the complement of a torus knot. Proc. Amer. Math. Soc. 15 (1964) 154--158.

\noindent
[K] R. Kirby, Problems in Low-Dimensional Topology.
Geometric topology.
Edited by William H. Kazez. AMS/IP Studies in Advanced
Mathematics, 2.2.  (1997) 35-473.

\noindent [Rs] R. Roberts, Constructing taut foliations. Comment.
Math. Helv. 70 (1995) 516--545.

\noindent [Rn] D. Rolfsen, Knots and Links. Publish or Perish 1976.

\noindent[V] J. Van Buskirk, Positive knots have positive Conway
polynomial.  Lecture Notes in Math. 1144 (1985) 146--159.

\noindent
[W] Y. Wu,
Dehn surgery on arborescent knots. J. Differential Geom. 43 (1996) 171--197.

\noindent[X]  P. Xu, The genus of closed 3-braids. J. Knot Theory
Ramifications 1 (1992) 303-326.

}
\end{document}